\documentclass{article}
\usepackage{amsmath,bm,colortbl,amssymb,url,amsthm}
\theoremstyle{plain}
\newtheorem{theorem}{Theorem}
\newcommand{\T}{\textsf{T}} 
\makeatletter 
\def\@biblabel#1{\hspace*{-\labelsep}}
\makeatother
\allowdisplaybreaks

\begin{document}
\title{The evolution of the game of baccarat}
\author{S. N. Ethier\thanks{Department of Mathematics, University of Utah, 155 South 1400 East, Salt Lake City, UT 84112, USA. e-mail: ethier@math.utah.edu.  Partially supported by a grant from the Simons Foundation (209632).}\ \ and 
Jiyeon Lee\thanks{Department of Statistics, Yeungnam University, 214-1 Daedong, Kyeongsan, Kyeongbuk 712-749, South Korea. e-mail: leejy@yu.ac.kr.  Supported by the Basic Science Research Program through the National Research Foundation of Korea (NRF) funded by the Ministry of Education (2013R1A1A3A04007670).}}
\date{}
\maketitle

\begin{abstract}
The game of baccarat has evolved from a parlor game played by French aristocrats in the first half of the 19th century to a casino game that generated over US\$41 billion in revenue for the casinos of Macau in 2013.  The parlor game was originally a three-person zero-sum game.  Later in the 19th century it was simplified to a two-person zero-sum game.  Early in the 20th century the parlor game became a casino game, no longer zero-sum.  In the mid 20th century, the strategic casino game became a nonstrategic game, with players competing against the house instead of against each other.  We argue that this evolution was motivated by both economic and game-theoretic considerations.
\end{abstract}

\section{Introduction}

The game of baccarat, which generated over US\$41 billion in revenue for the casinos of Macau in 2013, is a descendant of a 19th-century French parlor game.  But games evolve and a number of changes have been made to baccarat over the more than 170 years it is known to have been played.  Our aim in this paper is to explain the likely motivation for these changes.  In brief, the primary aim was to make the game more profitable for casino operators.  The identities of the people responsible for this evolution are, with one exception, lost to history, but there can be no doubt that they were remarkably successful in achieving their goal.

In Section \ref{history} we provide the historical background for the game of baccarat, emphasizing the three principal forms of the game, baccarat banque (a three-person game), baccarat chemin de fer (a two-person game), and baccarat punto banco (a nonstrategic game).  Section \ref{chemin-parlor} treats the zero-sum parlor game of baccarat chemin de fer in detail.  Section \ref{chemin-casino} considers the nonzero-sum casino game of baccarat chemin de fer, both in its original form and in its modern form.  Section \ref{punto} discusses baccarat punto banco, which is the form of baccarat most widely played today.  Finally, we summarize the reasons for game's evolution in Section \ref{summary}.

\section{Baccarat history}\label{history}

A number of authors have claimed that the game of baccarat (usually spelled ``baccara'' in French) dates back to 15th-century Italy, but this claim is questionable because the earliest known mention of the game appears in Van-Tenac (1847, pp.~84--97); on the other hand, the game was played exclusively in aristocratic circles, which may explain why it ``does not grace the realms of recorded history before the nineteenth century'' (Parlett 1991, pp.~81--82).

The Van-Tenac reference deals with baccarat banque, a three-person zero-sum parlor game, which is likely the earliest form of baccarat.  A simpler version of the parlor game, baccarat chemin de fer, is a two-person zero-sum game that appeared later.  The first reliable mathematical study is due to Dormoy (1873), who referred to the latter as baccarat tournant (rotating baccarat).  The term ``baccarat chemin de fer'' (baccarat railway) came into vogue by 1880 (Badoureau 1881).  Baccarat chemin de fer was studied by Bertrand (1889, pp.~38--42), who was unaware of Dormoy's (1873) monograph, and although Borel (1924) correctly characterized Bertrand's work as ``extremely incomplete,'' it motivated Borel to begin the development of game theory in the 1920s (Dimand and Dimand 1996, p.~132).  The first game-theoretic analysis of baccarat chemin de fer was carried out by Kemeny and Snell (1957).

Well before the legalization of casino gambling in France in 1907, baccarat was played in clubs, and to pay the expenses a commission was levied.  At baccarat chemin de fer, the amount of the commission may have varied to some extent (Villiod 1906, p.~165), but a five percent commission on Banker wins became standard.  This makes the game a two-person nonzero-sum game, or a bimatrix game, that has a unique Nash equilibrium, at least under the Kemeny--Snell ``with replacement'' assumption.    

Eventually, Banker's strategic options at baccarat chemin de fer were severely constrained, with the goal of ensuring optimal play, and perhaps reducing the level of skill needed to play the role of Banker.  Early descriptions of casino baccarat chemin de fer (e.g., Beresford 1926, p.~74) included the five percent commission on Banker wins, but Banker's strategy was still unconstrained.  It was also unconstrained in the descriptions of Le Myre (1935) and Boll (1936, 1944).  The constraints were present in Scarne's (1949, p.~206) description of the game, though he was referring not to French casinos but to (illegal) US casinos.  

In pre-revolutionary Havana, a third form of baccarat, sometimes called punto banco but more often just baccarat, was introduced.  In baccarat punto banco the casino banks the game, offering bets on Player and on Banker.  Moreover, the strategic options of baccarat chemin de fer are eliminated, resulting in a nonstrategic game, which has become the most widely played form of baccarat.  It currently enjoys its greatest popularity in the territory of Macau, where the 35 casinos generated gambling revenue of about US\$45 billion in 2013, of which over 91\% came from baccarat (Schwartz 2014).  Baccarat punto banco (the Spanish term is ``punto y banca'') is said to have first appeared in Argentina, but the Argentinian game was different:  Although Banker's strategy was fixed, Player still had the option of drawing or standing on 5 and the game was not house-banked (Chapitre 1954, pp.~25--26).  The Cuban innovations were to make the game house-banked, allowing bets on Player and on Banker, as well as to mandate Player drawing on 5.  It seems likely that the former modification occurred during the 1940s, and the latter during the 1950s.  Indeed, Scarne (1949, p.~213) referred to ``baccarat'' as a house-banked version of chemin de fer, though without mentioning where it was played.  And Tommy Renzoni wrote, ``I opened the Baccarat game at the Capri [in Havana, c.\ 1956] with a somewhat revised set of rules that made the game more automatic, more a matter of pure chance.'' (Renzoni and Knapp 1973, pp.~49--50.)  This suggests that it was he who mandated Player drawing on 5.  It is well known that this form of baccarat came to the Las Vegas Sands in November 1959 (Renzoni and Knapp 1973, p.~69), and to Macau in 1962.

Thus, over a lifetime of at least 170 years, the game of baccarat has undergone five significant rules changes (the three-person game of baccarat banque was simplified to the two-person game of baccarat chemin de fer;  a five percent commission was imposed on Banker wins at baccarat chemin de fer;  Banker's strategic options at baccarat chemin de fer were severely constrained; baccarat punto banco was introduced as a house-banked game, with bets offered on Player and on Banker; and all remaining strategic options were eliminated at baccarat punto banco), yet in most respects it is the same game.  Our aim here is to argue that the game's evolution was motivated by economic and game-theoretic considerations.  See Section \ref{summary} for a summary of our conclusions in this regard.

\section{The parlor game of baccarat chemin de fer}\label{chemin-parlor}

The rules of the parlor game of baccarat chemin de fer are as follows.  The role of Banker rotates among the players (counter-clockwise), changing hands after a Banker loss or when Banker chooses to relinquish his role.  Banker announces the amount he is willing to risk, and the total amount bet on Player's hand cannot exceed that amount.  After a Banker win, all winnings must be added to the bank unless Banker chooses to withdraw.  Denominations A, 2--9, 10, J, Q, K have values 1, 2--9, 0, 0, 0, 0, respectively, and suits are irrelevant.  The total of a hand, comprising two or three cards, is the sum of the values of the cards, modulo 10.  In other words, only the final digit of the sum is used to evaluate a hand.  Two cards are dealt face down to Player and two face down to Banker, and each looks only at his own hand.  The object of the game is to have the higher total (closer to 9) at the end of play.  Winning bets on Player's hand are paid by Banker at even odds.  Losing bets on Player's hand are collected by Banker.  Hands of equal total result in a tie or a \textit{push} (no money is exchanged).  A two-card total of 8 or 9 is a \textit{natural}.  If either hand is a natural, play ends.  If neither hand is a natural, Player then has the option of drawing a third card (but this option is heavily constrained; see below).  If he exercises this option, his third card is dealt face up.  Next, Banker, observing Player's third card, if any, has the option of drawing a third card.  This completes the game.  Since several players can bet on Player's hand, Player's strategy is restricted.  He must draw on a two-card total of 4 or less and stand on a two-card total of 6 or 7.  When his two-card total is 5, he is free to draw or stand as he chooses.  (The decision is usually made by the player with the largest bet.)  Banker, on whose hand no one can bet, has no constraints on his strategy.  

This is the original form of the game.  If we assume, as did Kemeny and Snell (1957), that cards are dealt with replacement from a single deck (often expressed by assuming a shoe with infinitely many decks), and that each of Player and Banker sees the total of his own two-card hand but not its composition, then baccarat is a $2\times2^{88}$ matrix game.  Let us reduce the size of the game considerably by eliminating strictly dominated columns of the payoff matrix $\bm A$.  We find that Banker's optimal move does not depend on Player's strategy in 84 of the 88 strategic situations, the exceptions being $(3,9)$ (Banker total 3, Player third card value 9), $(4,1)$, $(5,4)$, and $(6,\varnothing)$ (Banker total 6, Player stands).  See Table~\ref{reduction}.

\begin{table}[htb]
\caption{\label{reduction}Banker's optimal move at baccarat chemin de fer, indicated by D (draw) or S (stand), except in the four cases indicated by $*$ in which it depends on Player's strategy.  The shading of D entries is for improved readability.\medskip}
\catcode`@=\active\def@{\phantom{0}}
\catcode`#=\active\def#{\phantom{$^0$}}
\begin{center}
\begin{tabular}{ccccccccccccc}
\hline
\noalign{\smallskip}
Banker's &\multicolumn{11}{c}{Player's third card value ($\varnothing$ if Player stands)}\\
two-card &&&&&&&&&&\\
total & 0# & 1# & 2# & 3# & 4# & 5# & 6# & 7# & 8# & 9# & $\varnothing$ \\
\noalign{\smallskip} \hline
\noalign{\smallskip}
$0$& \cellcolor[gray]{0.85}D# & \cellcolor[gray]{0.85}D# & \cellcolor[gray]{0.85}D# & \cellcolor[gray]{0.85}D# & \cellcolor[gray]{0.85}D# & \cellcolor[gray]{0.85}D# & \cellcolor[gray]{0.85}D# & \cellcolor[gray]{0.85}D# & \cellcolor[gray]{0.85}D# & \cellcolor[gray]{0.85}D# & \cellcolor[gray]{0.85}D\\
$1$& \cellcolor[gray]{0.85}D# & \cellcolor[gray]{0.85}D# & \cellcolor[gray]{0.85}D# & \cellcolor[gray]{0.85}D# & \cellcolor[gray]{0.85}D# & \cellcolor[gray]{0.85}D# & \cellcolor[gray]{0.85}D# & \cellcolor[gray]{0.85}D# & \cellcolor[gray]{0.85}D# & \cellcolor[gray]{0.85}D# & \cellcolor[gray]{0.85}D\\
$2$& \cellcolor[gray]{0.85}D# & \cellcolor[gray]{0.85}D# & \cellcolor[gray]{0.85}D# & \cellcolor[gray]{0.85}D# & \cellcolor[gray]{0.85}D# & \cellcolor[gray]{0.85}D# & \cellcolor[gray]{0.85}D# & \cellcolor[gray]{0.85}D# & \cellcolor[gray]{0.85}D# & \cellcolor[gray]{0.85}D# & \cellcolor[gray]{0.85}D\\
3&       \cellcolor[gray]{0.85}D# & \cellcolor[gray]{0.85}D# & \cellcolor[gray]{0.85}D# & \cellcolor[gray]{0.85}D# & \cellcolor[gray]{0.85}D# & \cellcolor[gray]{0.85}D# & \cellcolor[gray]{0.85}D# & \cellcolor[gray]{0.85}D# & S# & $*$# & \cellcolor[gray]{0.85}D\\
4&       S# & $*$# & \cellcolor[gray]{0.85}D# & \cellcolor[gray]{0.85}D# & \cellcolor[gray]{0.85}D# & \cellcolor[gray]{0.85}D# & \cellcolor[gray]{0.85}D# & \cellcolor[gray]{0.85}D# & S# & S# & \cellcolor[gray]{0.85}D\\
5&       S# & S# & S# & S# & $*$# & \cellcolor[gray]{0.85}D# & \cellcolor[gray]{0.85}D# & \cellcolor[gray]{0.85}D# & S# & S# & \cellcolor[gray]{0.85}D\\
6&       S# & S# & S# & S# & S# & S# & \cellcolor[gray]{0.85}D# & \cellcolor[gray]{0.85}D# & S# & S# & $*$\\
7&       S# & S# & S# & S# & S# & S# & S# & S# & S# & S# & S\\
\noalign{\smallskip}
\hline
\end{tabular}
\end{center}
\end{table}

The case of $(3,9)$ occurred in the climactic hand in Ian Fleming's (1953) novel \textit{Casino Royale} (see Ethier 2010, pp.~616--617, for the full quotation), with James Bond in the role of Player, Le Chiffre in the role of Banker, and 32 million francs at stake.  As Fleming noted (Chap. 13), ``Holding a three and giving nine is one of the moot situations at the game. The odds are so nearly divided between to draw or not to draw.''  Le Chiffre chose to draw (which happens to be the correct move from the game-theoretic perspective), drawing 5 for a total of 8.  Then Bond's hand was turned over to reveal two queens for a total of 9 and the win.  Technically, the game they played was not quite baccarat chemin de fer because the role of Banker did not rotate among the players.  Instead, Le Chiffre purchased the role of Banker for one million francs (Chap.~9).  

Using Table~\ref{reduction} we can reduce the matrix game to $2\times2^4$, so we now regard $\bm A$ as being a $2\times16$ matrix.  We can further reduce this $2\times16$ matrix to a $2\times5$ matrix using strict dominance.  Specifically, $\bm A$ is replaced by the matrix (again denoted by $\bm A$)
\begin{equation}\label{A,2x5}
\bm A=\bordermatrix{& \text{SSSS} & \text{SSDS} & \text{DSDS} & \text{DSDD} & \text{DDDD}\cr
\text{S on 5} & -4636 & -4635 & -4564 & -2692 & -2585\cr              
\text{D on 5} & -3585 & -3600 & -3705 & -4121 & -4126\cr}\frac{16}{(13)^6},
\end{equation}
where, for example, the Banker pure strategy DSDS means draw on $(3,9)$, stand on $(4,1)$, draw on $(5,4)$, and stand on $(6,\varnothing)$.  It is easy to check that the kernel is specified by columns DSDS and DSDD, and that leads to the following theorem.

\begin{theorem}[Kemeny--Snell 1957]\label{KS-thm}
The parlor game of baccarat chemin de fer, a $2\times2^{88}$ matrix game, has a unique solution.  Player's optimal mixed strategy is to draw on $5$ with probability $p:=9/11$.  Banker's optimal mixed strategy is as in Table~\ref{reduction}, except that he draws on $(3,9)$, stands on $(4,1)$, draws on $(5,4)$, and mixes on $(6,\varnothing)$, drawing with probability $q:=859/2288$.  The value of the game (to Player) is $v_P:=-679568/[11(13)^6]\approx-0.0127991$.
\end{theorem}

See Ethier (2010, Chap.~5) for a detailed proof and Deloche and Oguer (2007b) for an alternative derivation.  Other analyses of the parlor game of baccarat chemin de fer have been given by Foster (1964), Downton and Lockwood (1975), and Ethier and G\'amez (2013).  These papers treated increasingly realistic models of the game.

\section{The casino game of baccarat chemin de fer}\label{chemin-casino}

Initially, the casino game differed from the parlor game in only one respect, the imposition of a five percent commission on Banker wins.  Let us consider, more generally, a 100$\alpha$ percent commission on Banker wins, where $0\le\alpha<1/15$.  If we continue to assume that cards are dealt with replacement from a single deck, and that each of Player and Banker sees the total of his own two-card hand but not its composition, then baccarat is a $2\times2^{88}$ bimatrix game, with payoff matrix $\bm A$ to Player as before, and payoff matrix $\bm B_\alpha$ to Banker with Banker winning $1-\alpha$ per unit lost by Player and losing 1 per unit won by Player.

Let us reduce the size of the game considerably by eliminating strictly dominated columns of $\bm B_\alpha$.  Assuming $0\le\alpha<1/15$, we find that Banker's optimal move does not depend on Player's strategy in 84 of the 88 strategic situations, the exceptions again being $(3,9)$, $(4,1)$, $(5,4)$, and $(6,\varnothing)$.  Table~\ref{reduction} still applies.  This reduces the game to a $2\times2^4$ bimatrix game.  We can further reduce it to $2\times5$, with $\bm A$ as in \eqref{A,2x5} and $\bm B_\alpha$ equal to the matrix whose transpose is
\begin{equation}\label{B,2x5}
\bm B_\alpha^{\T}=\bordermatrix{& \text{S on 5} & \text{D on 5}\cr
\text{SSSS} & 9272 - 278353 \alpha & 7170 - 276363 \alpha\cr 
\text{SSDS} & 9270 - 278423 \alpha & 7200 - 276433 \alpha\cr 
\text{DSDS} & 9128 - 278423 \alpha & 7410 - 276593 \alpha\cr 
\text{DSDD} & 5384 - 278007 \alpha & 8242 - 278673 \alpha\cr 
\text{DDDD} & 5170 - 277971 \alpha & 8252 - 278733 \alpha\cr}\frac{8}{(13)^6}.
\end{equation}

We know that there exists a Nash equilibrium, and we can use the support enumeration algorithm to find it and show it is unique.  This suffices, provided the game is \textit{nondegenerate}, which means that no mixed strategy of support size $s\ge1$ has more than $s$ pure best responses.  This hypothesis can easily be verified.

\begin{theorem}\label{Nash}
The classic casino game of baccarat chemin de fer with a $100\alpha$ percent commission on Banker wins, a $2\times2^{88}$ bimatrix game parameterized by $\alpha$, has a unique Nash equilibrium for $0\le\alpha<1/15$.  Player's equilibrium mixed strategy is to draw on $5$ with probability
\begin{equation}\label{p(alpha)}
p:=\frac{9-\alpha}{11-6\alpha}.
\end{equation}
Banker's equilibrium mixed strategy is as in Table~\ref{reduction}, except that he draws on $(3,9)$, stands on $(4,1)$, draws on $(5,4)$, and mixes on $(6,\varnothing)$, drawing with probability $q:=859/2288$.  The safety level for Player is $v_P$ of Theorem \ref{KS-thm} and that for Banker is
\begin{equation*}
v_B:=\frac{8 (84946 - 3099233 \alpha + 1668708 \alpha^2)}{(11 - 6 \alpha)(13)^6}.
\end{equation*}
If $\alpha=0.05$, then $v_B\approx-0.0101991$.  These are also the expected payoffs if both Player and Banker use their equilibrium strategies.
\end{theorem}

By \textit{safety level} we mean the amount a player can guarantee himself on average, regardless of his opponent's strategy.

Banker's safety level $v_B$ is greater than Player's safety level $v_P$ if and only if $0\le\alpha<\alpha'$, where 
$\alpha':=(34601239 - \sqrt{1060031672799697})/36711576\approx0.0556531$.  A commission on Banker wins greater than or equal to $\alpha'$ would be counter-productive because there would then be no incentive for players to take the role of Banker.  (We are assuming, in effect, that all players are rational and knowledgeable, which may be unrealistic.)  From the casino's perspective then, the commission $\alpha$ should be maximized subject to the constraints that $\alpha<\alpha'$ and $\alpha$ is simple enough to permit mental calculations.  5.5\% was presumably deemed too complicated, so five percent became the accepted figure.  

The reason for assuming $0\le\alpha<1/15$ in Theorem \ref{Nash} is that this is the maximal interval over which Table~\ref{reduction} applies without change.

As mentioned previously, in modern baccarat chemin de fer, Banker's strategy is highly constrained.  Specifically, the 84 of the 88 strategic situations that require a draw or stand decision in Table~\ref{reduction} are all part of the modern Banker strategy.  In addition, Banker stands on $(4,1)$ (Banker total 4, Player third card value 1), perhaps because the improvement in expected gain from drawing instead of standing when Player draws on 5 is very small (about 0.0025641 when $\alpha=0.05$).  Finally, Banker stands on $(6,\varnothing)$ (Banker total 6, Player stands), despite the fact that the improvement in expected gain from drawing instead of standing when Player draws on 5 is rather substantial (about 0.0673077 when $\alpha=0.05$).  Under these rules, which predate 1949 (Scarne 1949, p.~206) the only optional cases for Banker are $(3,9)$ and $(5,4)$.\footnote{We also find these rules in a 1920 Spanish book (Chapitre 1920), though they appear to be recommendations rather than requirements.  At $(3,9)$ and $(5,4)$ the recommendation is ``indifferent'', not ``optional''.}  Thus, we have a $2\times4$ bimatrix game with payoff bimatrix $(\bm A,\bm B_\alpha)$, where 
\begin{equation*}
\bm A=\bordermatrix{&\text{SS}&\text{SD}&\text{DS}&\text{DD}\cr
\text{S on 5} & -4636 & -4635 & -4565 & -4564 \cr
\text{D on 5} & -3585 & -3600 & -3690 & -3705 \cr}\frac{16}{(13)^6}
\end{equation*}
and
\begin{equation*}
\bm B_\alpha^{\T}=\bordermatrix{& \text{S on 5} & \text{D on 5}\cr
\text{SS} & 9272 - 278353 \alpha & 7170 - 276363 \alpha\cr 
\text{SD} & 9270 - 278423 \alpha & 7200 - 276433 \alpha\cr 
\text{DS} & 9130 - 278353 \alpha & 7380 - 276523 \alpha\cr 
\text{DD} & 9128 - 278423 \alpha & 7410 - 276593 \alpha\cr}\frac{8}{(13)^6};
\end{equation*}
here, for example, the Banker pure strategy DS means draw on $(3,9)$ and stand on $(5,4)$.

\begin{theorem}\label{Nash,2x4}
The modern casino game of baccarat chemin de fer with a $100\alpha$ percent commission on Banker wins, a $2\times4$ bimatrix game parameterized by $\alpha$, has a unique Nash equilibrium for $0\le\alpha<2/5$, and it is a pure Nash equilibrium.  Player's equilibrium strategy is to draw on $5$.  Banker's equilibrium strategy is as in Table~\ref{reduction}, except that he draws on $(3,9)$, stands on $(4,1)$, draws on $(5,4)$, and stands on $(6,\varnothing)$.  The safety levels are $v_P:=-3705(16)/(13)^6\approx-0.0122814$ and $v_B:=8(7410 - 276593\alpha)/(13)^6$.  If $\alpha=0.05$, then $v_B\approx-0.0106400$.  These are also the expected payoffs if both Player and Banker use their equilibrium strategies.
\end{theorem}

Because the casino levies a commission on Banker wins, it is to the casino's advantage to have Banker play well.  Perhaps this explains why Banker's strategy is so severely constrained in the modern game.  A more convincing explanation is that casinos may have wanted to reduce the level of skill required to take the role of Banker.  There has been speculation that Banker's strategy evolved gradually by trial and error, but that seems unlikely.  The essence of Table~\ref{reduction} dates back to the 19th century.  Specifically, Dormoy (1873), Billard (1883), and Hoffman (1891) had versions of Table~\ref{reduction}, but with errors in two or three entries.  Finally, by the 20th century, the exact values had been computed, for example by Le Myre (1935).  Thus, it is not surprising that, in the 84 entries that do not depend on Player's strategy, all appear correctly in the modern Banker strategy.

The only surprise is the stipulation that Banker stand on $(6,\varnothing)$ because this decision depends on Player's strategy.  (The same is true in the case of $(4,1)$, but we have noted that a mandated Banker stand makes sense in that case.)  Some authors have attributed this to a mistake, while others have suggested that it was an attempt to equalize the Player and Banker safety levels.  Interestingly, at Crockford's Club in London in the early 1960s, Banker was allowed to draw or stand on $(6,\varnothing)$ (Kendall and Murchland 1964) as well as on $(3,9)$ and $(5,4)$, and this was more than a decade after the mandatory stand on $(6,\varnothing)$ had become conventional wisdom.  Perhaps Crockford's had a mathematician advising them.

\section{The casino game of baccarat punto banco}\label{punto}

Although baccarat chemin de fer and baccarat banque are still available at the Casino de Monte-Carlo, the most widely played form of baccarat today is called baccarat punto banco (or just baccarat) and is not a strategic game because both Player and Banker have mandated strategies: Player draws on 5 or less and stands on 6 and 7, and Banker uses Table~\ref{reduction} and draws on $(3,9)$, stands on $(4,1)$, draws on $(5,4)$, and stands on $(6,\varnothing)$.  An even more important distinction between baccarat chemin de fer and baccarat punto banco is that, in the latter, one can bet either on Player (paid at even odds) or on Banker (paid at 19 to 20, equivalent to even odds with a five percent commission on winning Banker bets).  Both bets are banked by the casino.  

To better understand the flow of money at baccarat, let $P$ and $B$ be the probabilities that Player wins and Banker wins.  Then, under the usual ``with replacement'' assumption, $P=2153464/(13)^6$ and $B=2212744/(13)^6$, so
\begin{eqnarray}\label{P}
B-P&=&\frac{296400}{5(13)^6}\approx0.0122814,\\ \label{B}
P-(19/20)B&=&\frac{256786}{5(13)^6}\approx0.0106400,\\ \label{P+B}
(1/20)B&=&\frac{553186}{5(13)^6}\approx0.0229214.
\end{eqnarray}
Eq.~\eqref{P} is the casino's expected gain per unit bet on Player, \eqref{B} is the casino's expected gain per unit bet on Banker, and \eqref{P+B} is the casino's expected gain per unit bet at baccarat chemin de fer.  Notice that \eqref{P+B} is the sum of \eqref{P} and \eqref{B}.  This is consistent if we regard the casino's expected gain at baccarat chemin de fer as coming not just from Banker but from Player via Banker and from Banker directly.

In baccarat chemin de fer, the total amount bet on a hand is limited by the amount Banker is willing to risk; thus, the more bettors, the more likely that there will be unfulfilled demand.  To make this more explicit, let us suppose there are $n$ players who would be willing to bet $x_1,x_2,\ldots,x_n>0$ units on Player.  In baccarat chemin de fer, if Banker announces a bet of $y>0$ units, then the amount bet on Player (as well as the amount bet by Banker) would be $\min(s,y)$ units, where $s:=x_1+x_2+\cdots+x_n$;  whereas in baccarat punto banco there would be $s$ units bet on Player and $y$ units bet on Banker.  Thus, $s+y-2\min(s,y)=|s-y|$ units is the amount of unfulfilled demand.  Another advantage of having the casino bank the bets is that a game can begin with as few as one player; in baccarat chemin de fer, the game is delayed until the required number of players are available.

There is one attractive feature (from the casino's perspective) of baccarat chemin de fer:  The players may win or lose but the casino only wins; there is no risk to the house.  But as Lucas and Kilby (2012, p.~266) pointed out (again, from the casino's perspective), ``the feature that makes [baccarat chemin de fer] attractive is also its greatest weakness.''  Much more money is to be made by banking the game than by ``raking'' it, even if there is a resulting increase in risk.

There are at least two advantages to mandating the strategies of Player and Banker.  First, it avoids the unpleasantness of one player making a decision that negatively impacts other players.  Second, it speeds up the game when there are no decisions to be made.  The main issue we want to address here is why the mandatory drawing rules were chosen.  There are three potential explanations.

\begin{enumerate}
\item The Player and Banker strategies mandated by the rules of baccarat punto banco are the closest pair of pure strategies to the optimal mixed strategies of Kemeny and Snell (1957) for the parlor game of baccarat chemin de fer.

\item The Player and Banker strategies mandated by the rules of baccarat punto banco coincide with the pure Nash equilibrium at modern baccarat chemin de fer, in which Banker has options only on $(3,9)$ and $(5,4)$; see Section~\ref{chemin-casino}.  (When $\alpha=0$, the pure Nash equilibrium is a saddle point.)

\item The Banker strategy mandated by the rules of baccarat punto banco coincides with Banker's best response to Player's $(1/2,1/2)$ mixed strategy (that is, Player draws on 5 with probability 1/2), often expressed in the baccarat literature as Banker's best response when he does not know whether Player is in the habit of drawing or standing on 5.  Furthermore, the Player strategy mandated by the rules of baccarat punto banco coincides with Player's best response to Banker's mandated strategy.
\end{enumerate}

Since baccarat punto banco almost certainly predates the Kemeny--Snell paper, the first explanation is not credible and is merely a coincidence.  The second explanation is also doubtful because game-theoretic concepts such as saddle points were not well known in the casino industry in the 1950s and earlier.  Further, there is no mention of this property in the baccarat literature.  The third explanation is the most plausible.  Banker's best response to $(1/2,1/2)$ was implicitly noted in a well-known book by Le Myre (1935)\footnote{Although Le Myre's mixture was intended to be $(1/2,1/2)$ (p.~37), it was actually $(89/194,105/194)$ when Player draws (p.~85) and $(3/5,2/5)$ when Player stands (p.~51).} and explicitly noted by Boll (1936, \S 187) in another widely circulated book.  In a subsequent book, Boll (1944, Fig.~13) argued that a $(1/3,2/3)$ mixed strategy (that is, Player draws on 5 with probability 2/3) is more realistic, and this too has the same Banker best response.  This book was published the same year as von Neumann and Morgenstern's \textit{Theory of Games and Economic Behavior}, so it is not surprising that its author was unaware of the subject of game theory.  As a best response to a Player mixture of $(1/2,1/2)$, this Banker strategy had already been published, albeit slightly inaccurately, by Billard (1883) and Hoffmann (1891).  Finally, as is clear from Chapitre (1954, pp.~25--26), ``punto y banca'' had a mandated Banker strategy in Argentina, whereas Player still had a choice.  So when that game was modified for use in Havana, it would have been easy to check that Player's best response is to draw on 5.  This is consistent with Tommy Renzoni's statement mentioned earlier.  Thus, we can arrive logically at the mandated strategies at baccarat punto banco without having to assume that mistakes were made at some point in the distant past.

\section{Summary}\label{summary}

We mentioned that the game of baccarat has undergone five significant rules changes.  Here we review these changes and the likely reasons for them.  

The three-person game of baccarat banque was simplified to the two-person game of baccarat chemin de fer in the 19th century.  A likely reason is that baccarat banque required more players and was more difficult logistically.  Later, it would prove to be profitable in the very exclusive environment of high-stakes gambling (Graves 1963).

A five percent commission was imposed on Banker wins at baccarat chemin de fer to make it possible to offer the game in clubs and casinos.  The amount of the commission was presumably chosen as large as possible consistent with the goals of permitting simple mental calculations and not discouraging players from taking the role of Banker.  

Later, Banker's strategy was severely constrained, leaving only two draw-or-stand options, at $(3,9)$ and $(5,4)$, perhaps partly because it was in the interest of the casino to have Banker play well.  A more important reason is that less skill would be necessary to take the role of Banker.  The actual strategy chosen is very nearly optimal in the game-theoretic sense.  (Only the mandatory stands at $(4,1)$ and $(6,\varnothing)$ can be questioned, and as we have noted, at least the first of these has a logical justification.)  

Then, in Argentina, Banker's strategy became fixed, with mandatory draws at $(3,9)$ and $(5,4)$.  Perhaps the reason was to make baccarat more like blackjack, in which players do not have to fear a skillful dealer.  The game-theoretic justification for this fixed Banker strategy is that it is the best response to an equally weighted mixture of Player's two pure strategies, for which the concept (if not the implementation) dates back to the 19th century.  

The most important change to baccarat seems to have been made in Havana in the 1940s.  Baccarat punto banco became a house-banked game, offering bets on Player and on Banker.  This made the game vastly more profitable for the casinos (at the cost of additional risk), eliminating the need for offsetting bets, and allowing any number of players to bet as much as they want.  Further, the ability to bet on either Player or Banker gives the perception that the game is fair.

Finally, in baccarat punto banco as played in Havana in the 1950s, Player's strategy was also fixed by mandating a draw on 5, which is the best response to Banker's fixed strategy.  Eliminating strategy decisions speeds up the game and makes it purely a game of chance rather than one of chance and skill.  

All of these changes were motivated, at least in part, by a desire to make the game of baccarat more profitable for casino operators.  However, before a game can be profitable to casino operators it must be attractive to casino customers, and perceived fairness plays an important role in attractiveness.  Baccarat punto banco is almost unique among casino games in term of its perceived fairness and its profitability.

\section*{Acknowledgment}

We thank R\'egis Deloche for valuable advice concerning an earlier draft of this paper.

\end{document}